\documentclass{amsart}
\usepackage{amsmath}
\usepackage{amssymb}

  \usepackage{paralist}
  \usepackage{graphics} 
  \usepackage{epsfig} 
 \usepackage[colorlinks=true]{hyperref}
\hypersetup{urlcolor=blue, citecolor=red}

\newtheorem{theorem}{Theorem}[section]
\newtheorem{corollary}[theorem]{Corollary}

\newtheorem{lemma}[theorem]{Lemma}
\newtheorem{proposition}[theorem]{Proposition}
\newtheorem{conjecture}{Conjecture}
\newtheorem*{problem}{Problem}
\theoremstyle{definition}
\newtheorem{definition}[theorem]{Definition}
\newtheorem*{remark}{Remark}

\newcommand{\ep}{\varepsilon}

\newcommand{\dimb}{\mathrm{dim}_B}
\newcommand{\dimh}{\mathrm{dim}_H}
\newcommand{\calu}{\mathcal{U}}
\newcommand{\diam}{\mathrm{diam}}

\title[Decay of Lebesgue Numbers]
      {Exponential Decay of Expansive Constants}

\author[Peng Sun]{}


 \email{pengsunmath@hotmail.com}

\begin{document}

\maketitle

\centerline{\scshape Peng Sun }
\medskip
{\footnotesize
 \centerline{China Economics and Management Academy}
   \centerline{Central University of Finance and Economics}
   \centerline{No. 39 College South Road, Beijing, 100081, China}
} 

\bigskip

\begin{abstract}
A map $f$ on a compact metric space is expansive if and only if
$f^n$ is expansive. We study the exponential rate of decay of the expansive
 constant
and find some of its relations with other quantities about the dynamics,
such as dimension and entropy.
     
\end{abstract}

\section{Expansive maps}
Let $X$ be a compact metric space. A homeomorphism (continuous map)
 $f$ is called expansive
if there is $\gamma>0$ such that $d(f^n(x),f^n(y))<\gamma$ for all $n\in\mathbb{Z}$
($n\ge 0$) implies $x=y$. We call the largest $\gamma$ the expansive
constant of $f$, denoted by $\gamma(f)$, as it depends on $f$.

In this paper we assume that $f$ is an expansive homeomorphism, if not specified.
Results and their proofs for expansive continuous maps are very similar and
will be omitted.

Our discussion on expansive constants builds on the following proposition:

\begin{proposition}
For every $n\in\mathbb{Z}$ ($n\in\mathbb{N}$ for continuous maps), $f$ is
expansive if and only if $f^n$ is expansive.
\end{proposition}

\begin{proof}
From the definition, it is trivial that $f$ is expansive if and only if $f^{-1}$
 is expansive, and for $n\in\mathbb{N}$, if $f^n$ is expansive, then $f$
  is expansive.

Now assume that $f$ is expansive with expansive constant $\gamma(f)$. As
$f$ is continuous, there is $\epsilon>0$ such that $d(x,y)<\epsilon$ implies
$d_f^n(x,y)<\gamma(f)$, where
$$d_f^n(x,y)=\max_{0\le k\le n-1} d(f^k(x),f^k(y)).$$
So $d((f^n)^k(x),(f^n)^k(y))<\epsilon$ for all $k\ge 0$ implies
$d(f^{kn+j}(x),f^{kn+j}(y))<\gamma(f)$ for all $k\ge 0$ and $0\le j\le n-1$.
$f^n$ is expansive with expansive constant no less than $\epsilon$.
\end{proof}

From this proposition, if $f$ is expansive, we can talk about $\gamma(f^n)$,
the expansive constant of $f^n$. It is not difficult to make the following
observations:

\begin{lemma}
If $f$ is expansive, then the following hold.
\begin{enumerate}
\item (For homeomorphisms only) $\gamma(f)=\gamma(f^{-1})$.
\item For every $n\in\mathbb{N}$, $\gamma(f^n)\le\gamma(f)$.
\item If $d(x,y)<\tau$ implies $d_f^n(x,y)<\gamma(f)$, then $\gamma(f^n)\ge\tau$.
\end{enumerate}
\end{lemma}

It is natural to think about how $\gamma(f^n)$ varies as $n$ increases. For
Lipschitz maps we have a very rough estimate:

\begin{lemma}
If $f$ is expansive and Lipschitz with Lipschitz constant $L(f)$, then $L(f)>1$.
\end{lemma}

\begin{proposition}\label{lipmap}
Let $f$ be expansive and Lipschitz with Lipschitz constant $L(f)$, then for every $n\in\mathbb{N}$,
$\gamma(f^n)\ge\gamma(f)\cdot (L(f))^{-(n-1)}$.
\end{proposition}

\begin{proof}
$d(x,y)\le\gamma(f)\cdot (L(f))^{-(n-1)}$ implies $d_f^n(x,y)\le\gamma(f)$.
\end{proof}


\begin{corollary}\label{liphom}
If $f$ is a bi-Lipschitz homeomorphism with Lipschitz constants $L(f)$ and
$L(f^{-1})$, then for every $n\in\mathbb{N}$,
$$\gamma(f^n)\ge\gamma(f)\cdot\min_{0\le j\le n}\max\{(L(f))^{-j},(L(f^{-1}))^{-(n-j)}\}$$
\end{corollary}

\section{Exponential decay of expansive constants}
As in the usual way to study a quantity about the dynamics, we consider
the asymptotic exponential decay rate of the expansive constant $\gamma(f^n)$.

\begin{definition}
If $f$ is expansive, then let
$$h_E^+(f)=\limsup_{n\to\infty}-\frac1n\log\gamma(f^n)$$
$$h_E^-(f)=\liminf_{n\to\infty}-\frac1n\log\gamma(f^n)$$
\end{definition}

From now on we use $h_E^*$ to denote either $h_E^+$ or $h_E^-$, when the argument works for both cases. Some simple facts about them are listed below.

\begin{lemma}
If $f$ is expansive, then
\begin{enumerate}
\item $h_E^*(f)\ge 0$.
\item $h_E^-(f)\le h_E^+(f)$.
\item (For homeomorphisms only) $h_E^*(f)=h_E^*(f^{-1})$.
\end{enumerate}
\end{lemma}

\begin{proposition}
If $f$ is expansive, then for every $n\in\mathbb{N}$, $h_E^+(f^n)\le nh_E^+(f)$
 and $h_E^-(f^n)\ge nh_E^-(f^n)$.
\end{proposition}

\begin{proof}
$$h_E^+(f^n)=\limsup_{k\to\infty}-\frac1k\log\gamma(f^{kn})\le
n\limsup_{j\to\infty}-\frac1j\log\gamma(f^j)=nh_E^+(f).$$
The other one is analogous.
\end{proof}

\begin{proposition}\label{exlip}
If $f$ is continuous map that is expansive and Lipschitz, then $h_E^*(f)\le\log
L(f)$. If $f$ is homeomorphism that is expansive and bi-Lipschitz, then
$$h_E^*(f)\le\frac{\log L(f)\cdot\log L(f^{-1})}{\log L(f)+\log L(f^{-1})}.$$
In particular, if $L=\max\{L(f), L(f^{-1})\}$, then $h_E^*(f)\le\frac12\log
L$.
\end{proposition}

\begin{proof}
This is a corollary of Proposition \ref{lipmap} and Corollary \ref{liphom}.
\end{proof}




\begin{proposition}
Let $f:X\to X$ and $g:Y\to Y$ be expansive and there is a bi-Lipschitz conjugacy
$h:X\to Y$ between them. Then $h_E^*(f)=h_E^*(g)$.
\end{proposition}

\begin{proof}
For every $n$ and $k$, $d(g^{kn}(x), g^{kn}(y))<(L(h))^{-1}\cdot\gamma(f^n)$
implies 
$$d(h(f^{kn}(h^{-1}(x))), h(f^{kn}(h^{-1}(x))))<\gamma(f^n),$$ which
provides $h^{-1}(x)=h^{-1}(y)$, hence $x=y$. So $\gamma_n(g)\le (L(h))^{-1}\cdot\gamma(f^n)$.
Similar argument shows $\gamma(f^n)\le (L(h^{-1})^{-1}\cdot\gamma_n(g)$.
Taking limits, we have $h_E^*(f)=h_E^*(g)$.
\end{proof}

\begin{corollary}\label{eqmet}
$h_E^*(f)$ is invariant under strongly equivalent metrics.
\end{corollary}

\section{Relations with the exponential decay of Lebesgue numbers}
The exponential decay of Lebesgue numbers has been discussed in \cite{Su11a}.
We are somewhat surprised by the relation between it and the decay of expansive
constants we observe.

Recall that for every open cover $\calu$ of a compact metric space $X$, the
Lebesgue number $\delta(\calu)$ is defined as the largest positive number
such that every $\delta(\calu)$ ball is covered by some element of $\calu$.
Let $f$ be a continuous map on $X$. Let $\calu_f^n=\bigvee_{j=0}^{n-1}f^{-j}(\calu)$
and $\delta_n(f,\calu)=\delta(\calu_f^n)$. Define
$$h_L^-(f,\mathcal{U})=\liminf_{n\to\infty}-\frac1n\log\delta_n(f,\mathcal{U}),$$
$$h_L^+(f,\mathcal{U})=\limsup_{n\to\infty}-\frac1n\log\delta_n(f,\mathcal{U}),$$
$$h_L^-(f)=\sup_\calu h_L^-(f,\mathcal{U})$$
and
$$h_L^+(f)=\sup_\calu h_L^+(f,\mathcal{U}).$$ 

\begin{theorem}\label{HLE}
If $f$ is expansive, then $h_E^*(f)\le h_L^*(f)$.
\end{theorem}

\begin{proof}
Let $\calu$ be an open cover such that $\diam\calu<\gamma(f)$. Then for each
$n>0$, $d(x,y)<\delta_n(f,\calu)$ implies $d(f^j(x),f^j(y))<\diam\calu<\gamma(f)$
for $0\le j\le n-1$.
So if $d(f^{kn}(x),f^{kn}(y))<\delta_n(f,\calu)$ for every $k$,
then $d(f^m(x),f^m(y))<\gamma(f)$ for every $m=kn+j$, which runs over all
integers. This forces $x=y$ as $\gamma(f)$ is the expansive constant of $f$.
So $\gamma(f^n)\le\delta_n(f,\calu)$. Take the limit and we obtain
$h_E^*(f)\le h_L^*(f,\calu)\le h_L^*(f)$. (The last relation is in fact an
equality from \cite[Corollary 3.9]{Su11a})
\end{proof}

\section{Relations with entropy and dimension}
The most important fact we observe is that the product of $h_E^*(f)$ and
box dimension also bounds topological entropy. By Theorem \ref{HLE}, this
bound is (strictly, see Theorem \ref{hyptor} for example) better than 
\cite[Theorem 4.7]{Su11a}. But this result only makes sense when $f$ is
expansive.

\begin{theorem}\label{bdh}
Let $f$ be expansive on a compact metric space $X$. $\dimb^\pm X$ are the
upper and lower
box dimensions of $X$ and $h(f)$ is the topological entropy of $f$. Then
$h_E^-(f)\cdot\dimb^+ X\ge h(f)$ and $h_E^+(f)\cdot\dimb^-X\ge h(f)$.
\end{theorem}

\begin{proof}
We only show the first inequality. Proof of the other is similar.

The result is trivial if $h(f)=0$. Assume $h(f)>0$.
Take any $\lambda>\dimb^+ X$.
There is $\ep_0>0$ such that $\ep<\ep_0$ implies that there is an open cover
$\mathcal{U}$ of $X$ such that $\diam\mathcal{U}<\ep$ and $|\mathcal{U}|\le\ep^{-\lambda}$.
If $n$ large enough such that $\exp nh(f)=\exp h(f^n)>\ep^{-\lambda}$, then
$\mathcal{U}$ is not a generator under $f^n$ (see for example, \cite[Section
5.6]{PW}), as there are at most $\ep^{-k\lambda}$
elements in 
$\bigvee_{j=0}^{k-1}f^{-jn}(\calu)$, which makes
$$h(f^n,\calu)=\lim_{k\to\infty}\frac1k H(\bigvee_{j=0}^{k-1}f^{-jn}(\calu))\le\log(\ep^{-\lambda})<h(f^n).$$
There is $A\in\bigvee_{j=-\infty}^\infty f^{-jn}(\calu)$ that contains
at least two
points, say, $x$ and $y$. Then for every $j\in\mathbb{Z}$, $d(f^{jn}(x),f^{jn}(y))<\ep$.
$\ep$ is not an expansive constant for $f^n$ and $\gamma(f^n)<\ep$.

Now take $N>-\frac{\lambda\log\ep_0}{h(f)}$. For every $n>N$, take $\ep<\ep_0$
such that $\exp (n-1)h(f)<\ep^{-\lambda}<\exp nh(f)$. Then
$$-\frac{\log\gamma(f^n)}n>-\frac{\log\ep}n>\frac{(n-1)h(f)}{n\lambda}$$
Take the lower limit we get $h_E^-(f)\cdot\lambda>h(f)$. The result follows
since $\lambda$ is arbitrarily chosen.
\end{proof}

The idea of the proof is a byproduct of the following problem considered
by the author.

\begin{problem}
Let $f$ be an expansive homeomorphism on a compact metric space. What is
the smallest possible number of elements in a generator? How is this number
related to other properties of $f$?
\end{problem}

It is sure that this number should be at least $\exp h(f)$. So the best result
we can expect is the integer no less than $\exp h(f)$, and this is the case
for full shifts. However, even for subshifts of finite types we have
no idea at this moment.

As for subshifts of finite types, we observe the following fact.
\begin{proposition}
Let $\sigma_A$
be a subshift of finite type. Let $q>1$ and the metric on $\Omega_A$ be
defined by $d(\omega,\omega')=q^{-\min\{|j||\omega_j\ne\omega'_j\}}$.
Then $h_E^*(\sigma_A)\cdot\dimh\Omega_A=h(\sigma_A)$,
where $\dimh\Omega_A$ is the Hausdorff dimension of $\Omega_A$.
\end{proposition}

\begin{remark}
It is well-known that under the above assumptions, 
$$\dimb^+\Omega_A=\dimb^-\Omega_A=\dimh\Omega_A=\frac2{\log q}\lim_{k\to\infty}\frac1k\log\|A^k\|.$$
(or $\frac1{\log q}\lim_{k\to\infty}\frac1k\log\|A^k\|$ for one-sided shifts.)
\end{remark}

\begin{proof}
We only show the result for two-sided shifts. Proof for one-sided shifts
is analogous.

We know that $h(\sigma_A)=\lim_{k\to\infty}\frac1k\log\|A^k\|$. If $h(\sigma_A)=0$,
the result is trivial. Otherwise, it is enough to show that $h_E^*(\sigma_A)=\frac{\log
q}2$.

If $d(\omega,\omega')<q^{-n}$, then we must have $\omega_j=\omega'_j$ for
all $-n\le j\le n$. So $d(\sigma_A^{k(2n-1)}(\omega),\sigma_A^{k(2n-1)}(\omega))<q^{-n}$
for all $k\in\mathbb{N}$
implies $\omega=\omega'$, hence $\gamma(\sigma_A^{2n-1})\ge q^{-n}$. In fact,
$\gamma(\sigma_A^l)\ge q^{-n}$ for all $l\le 2n-1$. So
$h_E^-(\sigma_A)\ge\frac{\log q}2$.

By Theorem \ref{bdh}, $h_E^+(\sigma_A)\le\frac{\log q}2$. The result follows.
\end{proof}

The above result is not so difficult but very interesting. It seems that
 for certain expansive dynamical systems the topological entropy may also
  be given
by the product $h_E^*(f)\cdot\dimh(X)$. Though the box dimensions may be
different from the Hausdorff dimension when the metric is changed, we believe
that the result involving Hausdorff dimension is always true. Considering
Corollary \ref{eqmet}, the result is true for many other commonly used metrics
on symbolic spaces. A proof for all equivalent metrics is in process.

Moreover, as every Anosov diffeomorphism is expansive and has Markov partition,
we may expect the following result:

\begin{conjecture}
Let $f$ be an Anosov diffeomorphism on a $m$-dimensional compact manifold,
then $h(f)=m h_E^*(f)$.
\end{conjecture}

It would be a big surprise if the conjecture is true. Nevertheless, the following
fact might boost our confidence, at least a little bit.

\begin{theorem}\label{hyptor}
The conjecture is true for hyperbolic linear automorphisms on the 2-torus
with the standard metric.
\end{theorem}

\begin{proof}
Let $\lambda>1$ and $\lambda^{-1}$ be the eigenvalues. For every $\mu>\lambda$,
there is a metric that is strongly equivalent to the standard metric, such
that the diffeomorphism $f$ is a bi-Lipschitz map with Lipschitz constants
$L(f)<\mu$ and $L(f^{-1})<\mu$. By Proposition \ref{exlip} and Corollary
\ref{eqmet}, $h_E^*(f)<\frac12\log\mu$.
Since $\mu$ is arbitrarily taken and by Theorem \ref{bdh}, we have
$h_E^*(f)=\frac12\log\lambda=\frac12 h(f)$.
\end{proof}



\end{document}